\documentclass[a4paper,pdftex,reqno]{amsart}

\makeatletter
\@namedef{subjclassname@2020}{\textup{2020} Mathematics Subject Classification}
\makeatother

\usepackage{geometry}
\title[LM-construction and twisted Alexander invariants]{The Long--Moody construction and twisted Alexander invariants}
\author{Akihiro Takano}


\address{GRADUATE SCHOOL OF MATHEMATICAL SCIENCES, THE
UNIVERSITY OF TOKYO, 3-8-1 KOMABA, MEGURO-KU, TOKYO, 153-8914,
JAPAN}
\email{takano@ms.u-tokyo.ac.jp}

\usepackage{amsmath} 
\usepackage{amsfonts}
\usepackage{mathrsfs}
\usepackage{amssymb}
\usepackage{amsthm}
\usepackage{bm}
\usepackage[new]{old-arrows}
\usepackage{wrapfig}
\usepackage[dvipdfmx]{graphicx}
\usepackage{color}
\usepackage[labelsep=colon]{caption}
\usepackage[labelformat=empty,subrefformat=parens]{subcaption}
\usepackage{multicol}
\usepackage{units}
\usepackage[all]{xy}
\usepackage{tikz-cd}
\usepackage{hyperref}
\usepackage{mathpazo}

\newcommand{\relmiddle}[1]{\mathrel{}\middle#1\mathrel{}}

\newtheorem{thm}{Theorem}[section]

\newtheorem*{thm*}{Main Theorem}

\theoremstyle{definition}
\newtheorem{defi}[thm]{Definition}
\newtheorem{ex}[thm]{Example}
\newtheorem{remark}[thm]{Remark}

\newcommand{\bZ}{\mathbb{Z}}
\newcommand{\bR}{\mathbb{R}}
\newcommand{\bC}{\mathbb{C}}
\newcommand{\sB}{\mathscr{B}}
\newcommand{\sG}{\mathscr{G}}
\newcommand{\cT}{\mathcal{T}}
\newcommand{\LM}{\mathcal{LM}}
\newcommand{\Diag}{\mathrm{Diag}}
\newcommand{\cI}{\mathcal{I}}
\newcommand{\ft}{\mathfrak{t}}
\newcommand{\wtil}{\widetilde}
\newcommand{\dfr}{\displaystyle\frac}

\begin{document}

\begin{abstract}
In 1994, Long and Moody introduced a method to construct a new representation of the braid group from the representation of the braid group or the semidirect product of the braid group and the free group.
In this paper, we show that its matrix presentation is written using the Fox derivative, and also a relation with twisted Alexander invariants.
\end{abstract}

\maketitle

\section{Introduction} \label{intro}
When we are given a representation of the braid group $B_n$, it is natural to hope to construct a link invariant from the representation since any link is obtained as the closure of a braid.
For example, Burau \cite{Burau} constructed a representation $\sB \colon B_n \to GL_n (\bZ[t^{\pm1}])$ that reconstructs the Alexander polynomial.
Specifically, let $\wtil{\sB} \colon B_n \to GL_{n-1} (\bZ[t^{\pm1}])$ be the reduced Burau representation and $b \in B_n$ an $n$-braid, then
\begin{equation*}
\Delta_{\hat{b}} (t) = \pm t^l \frac{(t-1) \det(\wtil{\sB} (b) - I_{n-1})}{(t^n -1)},
\end{equation*}
where $l \in \bZ$, $\Delta_{\hat{b}} (t)$ is the Alexander polynomial of the closure of $b$, and $I_k$ is the $k \times k$ identity matrix.
Moreover, Birman \cite{Birman} showed that the multivariable Alexander polynomial $\Delta_{\hat{b}} (t_1, \ldots, t_n)$ of the closure of a pure braid $b$ is described by the reduced Gassner representation $\wtil{\sG} \colon P_n \to GL_{n-1} (\bZ[t_1^{\pm1}, \ldots, t_n^{\pm1}])$ of the pure braid group $P_n$ as follows:
\begin{equation*}
\Delta_{\hat{b}} (t_1, \ldots, t_n) = \pm t_1^{l_1} \cdots t_n^{l_n} \frac{\det(\wtil{\sG} (b) - I_{n-1})}{(t_1 \cdots t_n - 1)},
\end{equation*}
where $l_1, \ldots, l_n \in \bZ$.

As a generalization of the Alexander polynomial, Wada \cite{Wada} and Lin  \cite{Lin} independently defined the twisted Alexander invariant which is an invariant of a given link and a representation of the fundamental group of the link complement.
Similarly to the Alexander polynomial, there are various interpretations for this invariant, the Reidemeister-torsion \cite{Kitano}, the order of a certain twisted homology \cite{Kirk-Livingston}, and so on.
Conway \cite{Conway} introduced the twisted Burau map $\sB_{\rho} \colon B_c \to GL_{nk} (R[t_1^{\pm1}, \ldots, t_{\mu}^{\pm1}])$, where $R$ is a commutative ring and $\rho \colon F_n \to GL_k (R)$ is a representation of the free group $F_n$.
An $n$-braid $b$ is $\mu$-colored if each of its components is assigned an integer in $\{1, 2, \ldots, \mu \}$.
Then we obtain a sequence $c = (c_1, \ldots, c_n)$ of integers.
For such a sequence $c = (c_1, \ldots, c_n)$, we can define the subgroup $B_c$ of $B_n$, which is called the $c$-colored braid group; see Section \ref{colored}.
The twisted Burau map is defined as a homomorphism of the twisted homology of the $n$-punctured disk, and generally not a representation.
Conway showed an analogue of the above formulas, that is, the twisted Alexander invariant is obtained from the reduced twisted Burau map $\wtil{\sB}_{\rho} \colon B_c \to GL_{(n-1)k} (R[t_1^{\pm1}, \ldots, t_{\mu}^{\pm1}])$: suppose that $F_n$ is generated by $x_1, \ldots, x_n$ which are the standard generators of the fundamental group of the $n$-punctured disk.
For a $\mu$-colored braid $b \in B_c$, if a representation $\rho \colon F_n \to GL_k (R)$ factors through a representation $\pi_1 (S^3 \setminus \hat{b}) \to GL_k (R)$, then
\begin{equation*}
\Delta_{\hat{b}, \rho} (t_1, \ldots, t_{\mu}) = \pm \varepsilon t_{1}^{l_1} \cdots t_{\mu}^{l_{\mu}} \frac{\det(\wtil{\sB}_{\rho} (b) - I_{(n-1)k})}{ \det(\rho(x_1 \cdots x_n) t_{c_1} \cdots t_{c_n} -I_k)},
\end{equation*}
where $\varepsilon \in \det (\rho (\pi_1 (S^3 \setminus \hat{b})))$, $l_1, \ldots, l_{\mu} \in \bZ$, and $\Delta_{\hat{b}, \rho} (t)$ is the twisted Alexander invariant of the closure of $b$ associated with the representation $\pi_1 (S^3 \setminus \hat{b}) \to GL_k (R)$ induced by $\rho$.

On another note, Long and Moody \cite{Long} introduced a method of constructing a new representation of $B_n$ from a representation of $B_n \ltimes F_n$ or $B_{n+1}$, where the action of $B_n$ on $F_n$ is the Artin representation.
This method is called the Long--Moody construction.
More precisely, given a representation $\rho \colon B_n \ltimes F_n \to GL_{k}(R)$, its Long--Moody construction is an $nk$-dimensional representation of $B_n$.
This representation is more complicated and abundant than the initial one.
For example, if we take a one-dimensional trivial representation as the initial one, then we obtain the unreduced Burau representation.
Bigelow and Tian \cite{Bigelow-Tian} generalized this construction to some subgroups of $B_n$.
Souli\'e also studied this construction from a functorial point of view and extended it in \cite{Soulie1}, and then generalized it to other families of groups, such as the mapping class group of surfaces \cite{Soulie2} and the welded braid group \cite{Bellingeri-Soulie} with Bellingeri.

In this paper, we introduce the multivariable reduced Long--Moody construction $\wtil{\mathcal{LM}}_{\ft_c} (\rho) \colon B_c \to GL_{(n-1)k}(R[t_1^{\pm1}, \ldots, t_{\mu}^{\pm1} ])$ and show a relation with the twisted Alexander invariant.
In section \ref{twisted}, we introduce the definition of the twisted Alexander invariant with reference to Wada \cite{Wada}.
In section \ref{colored}, we recall the action of the braid group on the free group and explain the colored braids.
In section \ref{LM}, we define the multivariable Long--Moody construction which is a generalization of \cite{Bigelow-Tian} and show that a matrix presentation of the Long--Moody construction is written by the Fox derivative.
Then we give a formula related to the twisted Alexander invariant:
\setcounter{section}{4}
\setcounter{thm}{7}
\begin{thm*}[\textbf{Theorem \ref{main}}]
Let $b \in B_c$ be a $\mu$-colored braid and $\hat{b}$ its closure.
Let $\rho \colon B_c \ltimes F_n \to GL_{k}(R)$ be a representation such that the restriction $\rho|_{F_n} \colon F_n \to GL_{k}(R)$ factors through $\pi_1 (S^3 \setminus \hat{b}) \to GL_k (R)$.
Then
\begin{equation*}
\Delta_{\hat{b},\rho}(t_1, \ldots, t_{\mu}) = \pm \varepsilon t_{1}^{l_1} \cdots t_{\mu}^{l_{\mu}} \frac{\det \left( \wtil{\LM}_{\ft_c}(\rho) (b) - \Diag \left( \rho (b), \ldots, \rho (b) \right) \right)}{\det(\rho(x_1 \cdots x_n) t_{c_1} \cdots t_{c_n} -I_k)}
\end{equation*}
for some $\varepsilon \in R^{\times}$ and $l_1, \ldots, l_{\mu} \in \bZ$.
\end{thm*}
\setcounter{section}{1}
In section \ref{ex}, we give examples of the main theorem in the case of some knots: the trefoil knot, $(2,q)$-torus knot, and figure eight knot.

\section{Twisted Alexander invariants} \label{twisted}
Let $G$ be a group with a finite presentation
\begin{equation*}
G=\langle x_1 , \ldots , x_n \mid r_1 , \ldots , r_m \rangle
\end{equation*}
such that there exists a surjective homomorphism
\begin{equation*}
\alpha \colon G \to H := \bZ^{\oplus \mu} \cong  \langle t_1, \ldots, t_{\mu} \mid t_i t_j = t_j t_i\ (1 \leq i < j \leq \mu) \rangle.
\end{equation*}
Suppose that $G$ is the fundamental group of a $\mu$-component link complement, then we fix a homomorphism $\alpha \colon G \to H$ that maps the meridian of $i$-th component to $t_i$.
Let $\rho \colon G \to GL_k(R)$ be a representation, where $R$ is a unique factorization domain.
These maps are naturally extended to the ring homomorphisms $\wtil{\alpha} \colon \bZ[G] \to \bZ[H]$ and $\wtil{\rho} \colon \bZ[G] \to M_k(R)$, where $M_k(R)$ is the matrix algebra of degree $k$ over $R$.
Then the tensor product homomorphism $\wtil{\rho} \otimes \wtil{\alpha} \colon \bZ[G] \to M_k(R[H])$ of $\wtil{\rho}$ and $\wtil{\alpha}$ is defined by
\begin{equation*}
(\wtil{\rho} \otimes \wtil{\alpha}) (g) := \rho (g) \alpha (g)
\end{equation*}
for any $g \in G$.
Let $F_n$ be the free group generated by $x_1, \ldots, x_n$ and $\phi \colon F_n \to G$ the surjective homomorphism induced by the presentations.
Similarly to $\alpha$ and $\rho$, the homomorphism $\phi$ induces a ring homomorphism $\wtil{\phi} \colon \bZ[F_n] \to \bZ[G]$.
Then we obtain the ring homomorphism
\begin{equation*}
\Phi := (\wtil{\rho} \otimes \wtil{\alpha}) \circ \wtil{\phi} \colon \bZ[F_n] \to M_k(R[H]).
\end{equation*}

We define the $m \times n$ matrix $M$ whose $(i,j)$ component is the $k \times k$ matrix
\begin{equation*}
\Phi \left( \frac{\partial r_i}{\partial x_j} \right) \in M_k(R[H]),
\end{equation*}
where $\dfr{\partial}{\partial x_j}$ is the Fox derivative with respect to $x_j$, that is, a linear map $\bZ[F_n] \to \bZ[F_n]$ over $\bZ$ satisfying two conditions:
\begin{itemize}
\item $\dfr{\partial x_i}{\partial x_j} = \delta_{ij}$, where $\delta_{ij}$ is the Kronecker delta, and
\item $\dfr{\partial (g g')}{\partial x_j} = \dfr{\partial g}{\partial x_j} + g \dfr{\partial g'}{\partial x_j}$ for any $g, g' \in F_n$.
\end{itemize}
This matrix $M$ is called the \textbf{Alexander matrix} \cite[Page 244]{Wada} of the presentation of $G$ associated with the representation $\rho$.

For $1 \leq j \leq n$, let $M_j$ be the $m \times (n-1)$ matrix obtained from $M$ by removing the $j$-th column. We regard $M_j$ as an $mk \times (n-1)k$ matrix with entries in $R[H]$. For an $(n-1)k$-tuple of indices
\begin{equation*}
I=\left( i_1 , \ldots , i_{(n-1)k} \right)\ \ \left(1 \leq i_1 < \cdots < i_{(n-1)k} \leq mk \right),
\end{equation*}
we write $M^I_j$ for the $(n-1)k \times (n-1)k$ matrix consisting of the $i_l$-th rows of the matrix $M_j$, where $1 \leq l \leq (n-1)k$.

Since the map $\alpha \colon G \to H$ is surjective, there exists an integer $1 \leq i \leq n$ such that $\alpha(x_i)$ is non-trivial.
Therefore we have $\det \Phi (x_i - 1) \neq 0$ \cite[Lemma 2]{Wada}.
Moreover, for such integers $i, j$ and any choice of the indices $I$, the following equality holds \cite[Lemma 3]{Wada}:
\begin{equation*}
\det M^I_i \det \Phi (x_j - 1) = \pm \det M^I_j \det \Phi (x_i - 1).
\end{equation*}

\begin{defi}[{\cite[Page 246]{Wada}}]
The \textbf{twisted Alexander invariant} $\Delta_{G, \rho} (t_1, \ldots, t_{\mu})$ of the group $G$ associated with the representation $\rho$ is defined to be a rational expression
\begin{equation*}
\Delta_{G, \rho} (t_1, \ldots, t_{\mu}) := \frac{\gcd_I (\det M^I_j)}{\det \Phi (x_j - 1)}
\end{equation*}
provided that $\det \Phi (x_j - 1) \neq 0$.
This is well-defined up to a factor $\varepsilon t_1^{\l_1} \cdots t_{\mu}^{l_{\mu}}$, where $\varepsilon \in R^\times$ and $l_1, \ldots, l_{\mu} \in \mathbb{Z}$.
If $m < n - 1$, we define $\Delta_{G, \rho}(t_1, \ldots, t_{\mu}) := 0$.
\end{defi}

The twisted Alexander invariant is an invariant of the group $G$, the surjective homomorphism $\alpha$, and the representation $\rho$.
More precisely, the twisted Alexander invariant is independent of the choice of the presentation of $G$.
Furthermore, if two representations $\rho_1$ and $\rho_2$ are equivalent, then we have $\Delta_{G, \rho_1}(t_1, \ldots, t_{\mu}) \equiv \Delta_{G, \rho_2}(t_1, \ldots, t_{\mu})$.

\section{Colored braids} \label{colored}
In this section, we make recollections on braid groups and their natural actions on free groups.
We refer the reader to \cite{Kassel-Turaev} for further details.
Also we explain the notion of colored braid.

\begin{figure}[tbp]
\begin{center}
\includegraphics[height=90pt]{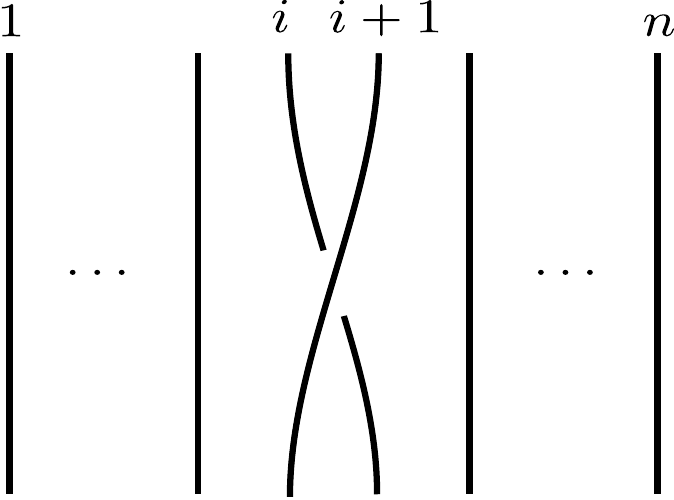}
\caption{The generator $\sigma_i$}
\label{braid_diagram}
\end{center}
\end{figure}

Let $B_n$ be the braid group of $n$ strands.
This group has the following presentation, which is called the Artin presentation:
\begin{equation*}
\left\langle \sigma_1 , \ldots , \sigma_{n-1} \relmiddle|
\begin{array}{ll}
\sigma_i \sigma_j = \sigma_j \sigma_i & (|i-j| \geq 2) \\ 
\sigma_i \sigma_{i+1} \sigma_i = \sigma_{i+1} \sigma_i \sigma_{i+1} & (1 \leq i \leq n-2)
\end{array}
\right\rangle.
\end{equation*}
The generator $\sigma_i$ corresponds to the $n$-braid described in Figure \ref{braid_diagram}.
There is a natural surjection from $B_n$ onto the symmetric group $S_n$.
Its kernel $P_n$ is called the \textbf{pure braid group}.

\begin{figure}[tbp]
\begin{center}
\includegraphics[height=140pt]{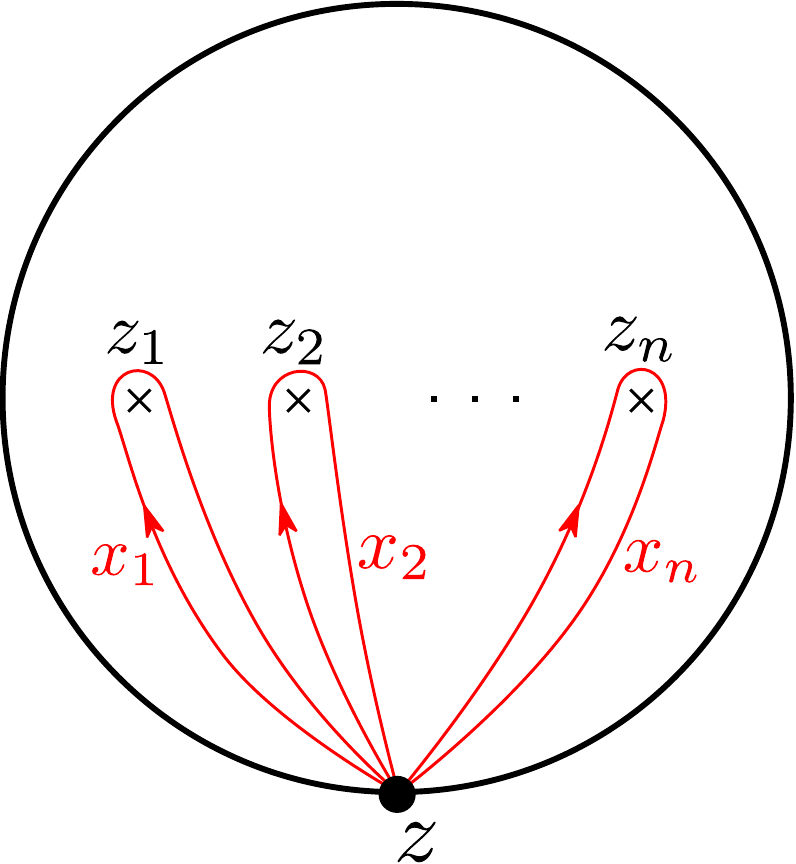}
\caption{The $n$-punctured disk $D_{n}$ and the generators of its fundamental group}
\label{disk_generators}
\end{center}
\end{figure}

\begin{figure}[tbp]
\begin{center}
\includegraphics[height=110pt]{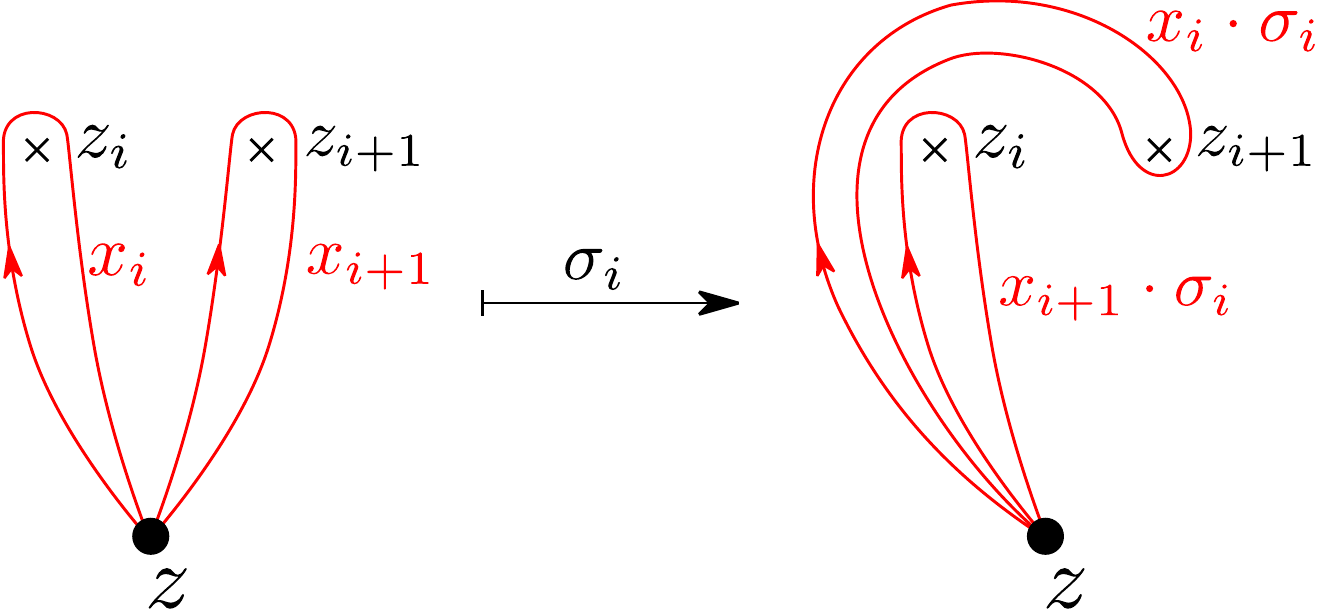}
\caption{The action of $\sigma_i$ on $\pi_1 (D_n, z)$}
\label{braid_action_disk}
\end{center}
\end{figure}

Let $D^2$ be the unit 2-disk in $\bR^2$ and fix $n$ distinct points $z_1, \ldots, z_n$ in the interior of $D^2$.
We shall assume so that each $z_i$ lies in $\textrm{Int} (D^2) \cap (\bR \times \{ 0 \}) = (-1,1) \times \{ 0 \}$ and $z_1 < \cdots < z_n$.
Set $D_n := D^2 \setminus \{ z_1, \ldots, z_n \}$ and fix a base point $z \in \partial D_n$.
Let $x_i$ be the simple loop around $z_i$ based at $z$ for $1 \leq i \leq n$ in Figure \ref{disk_generators}.
Then the fundamental group $\pi_1 (D_n, z)$ is identified with the free group $F_n = \left\langle x_1, \ldots, x_n \right\rangle$.
The braid group $B_n$ is isomorphic to the mapping class group $\mathcal{M} (D_n)$ of $D_n$, which is the group of isotopy classes of orientation-preserving self-homeomorphisms of $D_n$ fixing the boundary pointwise.
Note that the loop $x_1 \cdots x_n$ is homotopic to the boundary of $D_n$.
Since each element of $\mathcal{M} (D_n) \cong B_n$ induces the automorphism of $\pi_1 (D_n, z) \cong F_n$ as described in Figure \ref{braid_action_disk}, we have a \textit{right} action of $B_n$ on $F_n$ defined by
\begin{equation*}
x_j \cdot \sigma_i :=
\begin{cases}
x_i x_{i+1} x_i^{-1} & (j = i) \\
x_i & (j = i+1) \\
x_j & (j \neq i,i+1)
\end{cases}.
\end{equation*}

\begin{figure}[tbp]
\begin{center}
\includegraphics[height=120pt]{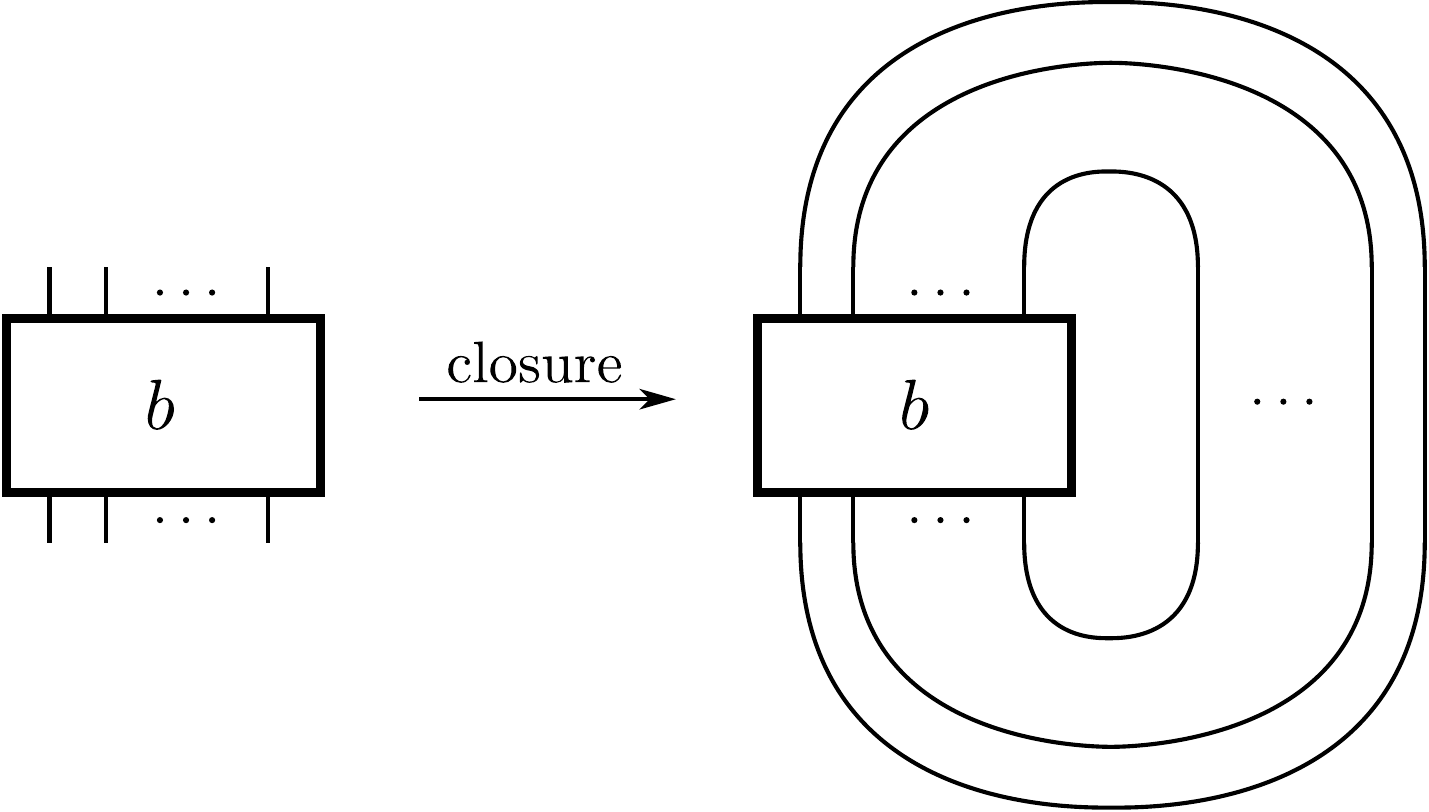}
\caption{The closure of a braid}
\label{braid_closure}
\end{center}
\end{figure}

For a braid $b \in B_n$, its \textbf{closure} $\hat{b}$ is the link in $S^3$ as described in Figure \ref{braid_closure}.
Then the fundamental group $\pi_1 (S^3 \setminus \hat{b})$ of its complement admits a presentation as follows:
\begin{equation*}
\left\langle x_1 , \ldots , x_n \relmiddle| x_1 = x_1 \cdot b, \ldots, x_n = x_n \cdot b \right\rangle.
\end{equation*}

\begin{figure}[tbp]
\begin{center}
\includegraphics[height=120pt]{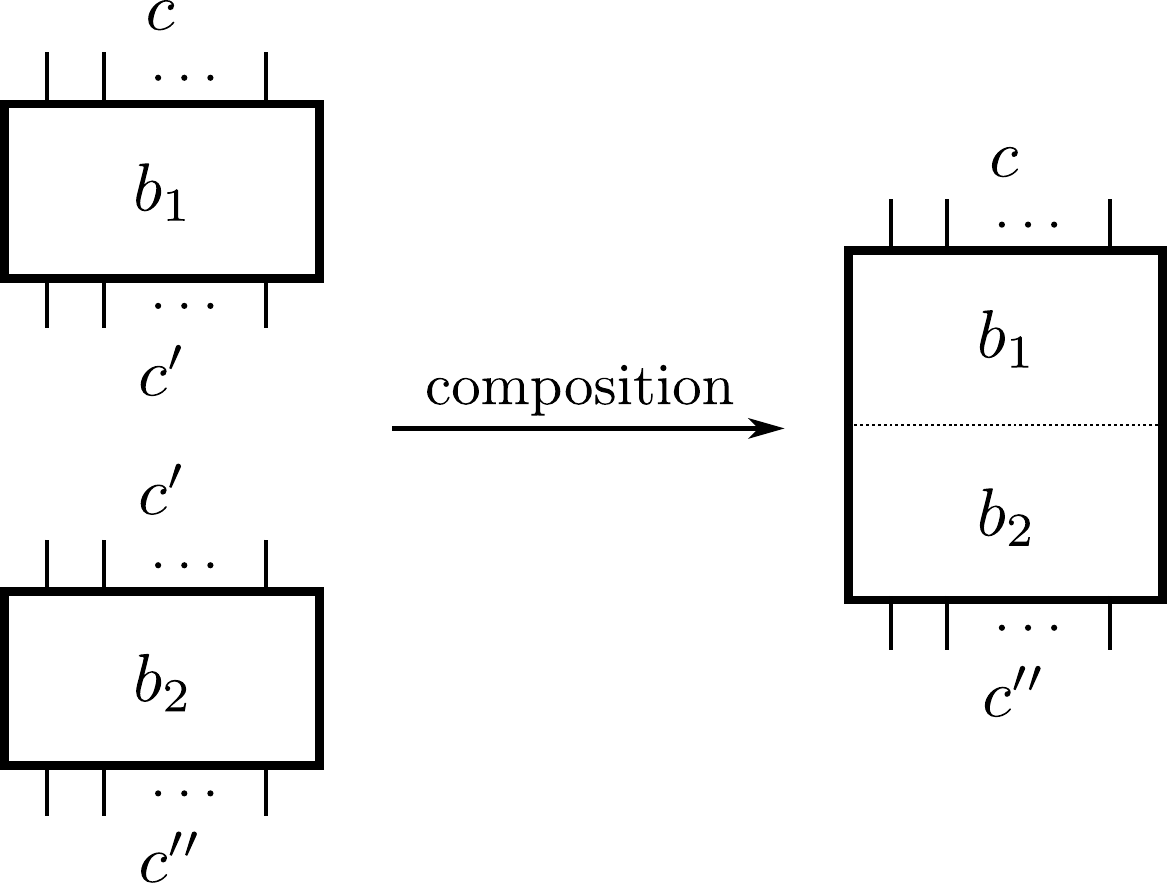}
\caption{$b_1$ is a $(c, c')$-braid, $b_2$ is a $(c', c'')$-braid, and their composition $b_1 b_2$ is a $(c, c'')$-braid.}
\label{colored_composition}
\end{center}
\end{figure}

An $n$-braid $b$ is \textbf{$\mu$-colored} if each of its components is assigned an integer in $\{1, 2, \ldots, \mu \}$ by a map $\{ 1, 2, \ldots, n \} \to \{1, 2, \ldots, \mu \}$.
Then we obtain a sequence $c = (c_1, \ldots, c_n)$ of integers, that is, $c_i \in \{1,2,  \ldots, \mu \}\ (1 \leq i \leq n)$ is the color of the $i$-th string.
The color of a braid induces a coloring on the $n$ initial points and $n$ terminal points.
The initial points are colored by the same as the sequence $c$, and let $c'$ be the coloring on the terminal points.
A $\mu$-colored braid is also called a $(c, c')$-braid.
If $b_1$ is a $(c, c')$-braid and $b_2$ a $(c', c'')$-braid, then their composition $b_1 b_2$ is a $(c, c'')$-braid; see Figure \ref{colored_composition}.
For any sequence $c$, the set $B_c$ is defined as the set of all $(c, c)$-braids.
This set has a group structure and is a subgroup of $B_n$.
The group $B_c$ is called the \textbf{$c$-colored braid group}.
For example, $B_{(1, 1, \ldots, 1)}$ is equal to $B_n$ and $B_{(1, 2, \ldots, n)}$ is equal to $P_n$.

\section{The Long--Moody construction} \label{LM}
In this section, we introduce the multivariable Long--Moody construction for the $c$-colored braid group $B_c$ to describe a relation with the (multivariable) twisted Alexander invariant.
This is a generalization of the one discussed in \cite{Bigelow-Tian}.

\subsection{\textit{Definition}}
We fix a sequence $c = (c_1, \ldots, c_n)$ of integers $1, 2, \ldots, \mu$, and set $H' := \langle t \rangle \oplus H \cong \bZ^{\oplus \mu+1}$.
We define a homomorphism $\ft_c \colon B_c \ltimes F_n \to H'$ by the following:
\begin{equation*}
\begin{cases}
b = \sigma_{i_1}^{m_1} \cdots \sigma_{i_l}^{m_l} \mapsto t^{m_1 + \cdots + m_l} & (b \in B_c)\\
x_j \mapsto t_{c_j} & (j = 1, \ldots, n)
\end{cases}.
\end{equation*}
From \cite[Theorem 3.1]{Bigelow-Tian}, we are able to apply the Long--Moody construction for a representation of $B_c \ltimes F_n$ that is a subgroup of $B_n \ltimes F_n$.
Let $\rho \colon B_c \ltimes F_n \to GL_k(R)$ be a representation, where $R$ is a commutative ring.
We regard $v \in R^{\oplus k}$ as a row vector $v = (v_1, \ldots, v_{k})$, and thus the representation $\rho$ is multiplied from the right of $v$.
Restricting $\rho$ and $\ft_c$ to the free group $F_n$ gives $R[H']^{\oplus k}$ a right $R[F_n]$-module structure.
Then we consider a $R[H']$-module $R[H']^{\oplus k} \otimes_{R[F_n]} \mathcal{I}_{F_n}$, where $\cI_{F_n}$ is the augmentation ideal of $F_n$, that is, $\cI_{F_n}$ is the kernel of the augmentation $R[F_n] \to R$.

\begin{defi}[cf. {\cite[Theorem 2.1]{Long}, \cite[Theorem 3,1]{Bigelow-Tian}}]
The \textbf{multivariable Long--Moody construction} of  the representation $\rho \colon B_c \ltimes F_n \to GL_k(R)$ is a representation $\LM_{\ft_c}(\rho) \colon B_c \to GL_{R[H']}(R[H']^{\oplus k} \otimes_{R[F_n]} \cI_{F_n})$ given by
\begin{equation*}
\LM_{\ft_c} (\rho) (b) (v \otimes h) := \bar{\ft}_c (b) \cdot \left( v \left((\rho \otimes \ft_c) (b) \right) \otimes h \cdot b \right)
\end{equation*}
for any $b \in B_c$, $h \in \cI_{F_n}$ and $v \in R[H']^{\oplus k}$, and the right action of $B_c$ on $F_n$ is naturally extended to the action on $\cI_{F_n}$.
Also, a homomorpshim $\bar{\ft}_c \colon B_c \to H'$ is defined as $\bar{\ft}_c (b) = \ft_c(b)^{-1}$ for any $b \in B_c$.
\end{defi}

By using the Fox derivative, we are able to calculate a matrix presentation for the Long--Moody construction:

\begin{thm} \label{matrix}
For any $b \in B_c$,
\begin{equation*}
\LM_{\ft_c}(\rho) (b) = \Diag \left( \rho (b), \ldots, \rho (b) \right) \cdot \left( (\rho \otimes \ft_c) \left( \frac{\partial (x_i \cdot b)}{\partial x_j} \right) \right)_{1 \leq i,j \leq n}.
\end{equation*}
\end{thm}

\begin{proof}
Note that $\cI_{F_n}$ is isomorphic to the free $R[F_n]$-module of rank $n$ generated by $x_i-1$ $(1 \leq i \leq n)$.
Therefore, there is an isomorphism
\begin{align*}
\begin{array}{ccc}
R[H']^{\oplus k} \otimes_{R[F_n]} \cI_{F_n} & \to & \displaystyle{\bigoplus_{i=1}^{n}} R[H']^{\oplus k} \vspace{5pt} \\
v \otimes (x_i-1) & \mapsto & 0 \oplus \cdots \oplus \underbrace{v}_{i-th} \oplus \cdots \oplus 0,
\end{array}
\end{align*}
and $GL_{R[H']}(R[H']^{\oplus k} \otimes_{R[F_n]} \cI_{F_n}) \cong GL_{nk}(R[H'])$.
Hence, for any $v_1, \ldots, v_n \in R[H']^{\oplus k}$, we identify $v_1 \oplus \cdots \oplus v_n \in \displaystyle{\bigoplus_{i=1}^{n}} R[H']^{\oplus k}$ with
\begin{equation*}
(v_1 \otimes (x_1-1)) + \cdots + (v_n \otimes (x_n-1)) \in R[H']^{\oplus k} \otimes_{R[F_n]} \cI_{F_n}.
\end{equation*}
By the fundamental formula of the Fox derivative,
\begin{equation*}
w - 1 = \sum^n_{j=1} \frac{\partial w}{\partial x_j} (x_j - 1)
\end{equation*}
for any $w \in \bZ [F_n]$.
Hence
\begin{align*}
\LM_{\ft_c}(\rho) (b) (v \otimes (x_i - 1)) &= \bar{\ft}_c (b) \cdot \left( v \left((\rho \otimes \ft_c) (b) \right) \otimes (x_i - 1) \cdot b \right) \\
&= \ft_c (b)^{-1} \cdot \ft_c (b) \cdot v (\rho (b)) \otimes \sum^n_{j=1} \left( \frac{\partial (x_i \cdot b)}{\partial x_j} \right) (x_j - 1) \\
&= v (\rho (b)) \otimes \sum^n_{j=1} \left( \frac{\partial (x_i \cdot b)}{\partial x_j} \right) (x_j - 1) \\
&= \sum^n_{j=1} v \left( \rho (b) \right) \left( (\rho \otimes \ft_c) \left( \frac{\partial (x_i \cdot b)}{\partial x_j} \right) \right) \otimes (x_j - 1).
\end{align*}
Therefore, we obtain the desired matrix presentation.
\end{proof}

If we consider the map $B_n \ni b \mapsto v \left((\rho \otimes \ft_c) (b) \right) \otimes h \cdot b$ for $h \in \cI_{F_n}$ and $v \in R[H']^{\oplus k}$, then this is also well-defined.
However, the variable $t$ remains in the image, and this is not essential.
Therefore, we take the product with the ``inverse'' homomorphism $\bar{\ft}_c$ to delete the variable $t$ in the image of $\LM_{\ft_c}(\rho)$.
Then $\LM_{\ft_c}(\rho)$ can be regarded as a representation on $GL_{R[H]} (R[H]^{\oplus k} \otimes_{R[F_n]} \cI_{F_n}) \cong GL_{nk}(R[H])$.
If $c = (1, 1, \ldots, 1)$, then this is a representation of $B_n$ which is equivalent to the representation of \cite[Corollary 2.6]{Long}.
Also, the case of $c = (1, 2, \ldots, n)$ is just as a representation of $P_n$ of \cite[Theorem 4.1]{Bigelow-Tian}.

\begin{ex}
Consider the one-dimensional trivial representation $\cT \colon B_c \ltimes F_n \to GL_1(\bZ)$.
Then the multivariable Long--Moody construction $\LM_{\ft_c}(\cT) \colon B_c \to GL_n(\bZ[t_1^{\pm1}, \ldots, t_{\mu}^{\pm1}])$ can be defined.
If $c = (1, 1, \ldots, 1)$, this is equivalent to the unreduced Burau representation.
Also, if $c = (1, 2, \ldots, n)$, this is equivalent to the unreduced Gassner representation (cf. \cite[Section 4]{Bigelow-Tian}).
\end{ex}

\begin{remark}
The original Long--Moody construction \cite[Theorem 2.1]{Long} is defined for $B_n$ without the homomorphism $\ft_{c}$.
In other words, the \textbf{Long--Moody construction} of  the representation $\rho \colon B_n \ltimes F_n \to GL_k(R)$ is a representation $\LM(\rho) \colon B_n \to GL_R(R^{\oplus k} \otimes_{R[F_n]} \cI_{F_n}) \cong GL_{nk} (R)$ given by
\begin{equation*}
\LM(\rho) (b) (v \otimes h) := v(\rho (b)) \otimes h \cdot b.
\end{equation*}
Moreover, this construction is also defined for the braid group $B_{n+1}$ via two ``canonical'' homomorphisms $F_n \to B_{n+1}$ and $ B_n \to Aut(F_n)$; see \cite[Theorem 2.4]{Long}.
Souli\'e \cite{Soulie1, Soulie2} extended and generalized this method.
For example, if we choose homomorphisms $\chi \colon F_n \to B_{n+1}$ and $\kappa \colon B_n \to Aut(F_n)$ satisfying the equation
\begin{equation*}
\chi (x_j) \sigma_{i+1} = \sigma_{i+1} \chi (\kappa(\sigma_i) (x_j)),
\end{equation*}
then the Long--Moody construction can be defined.
Since this construction depends on $\chi$ and $\kappa$, we may write it as $\LM_{\kappa, \chi} \colon B_n \to GL_R(R^{\oplus k} \otimes_{R[F_n]} \cI_{F_n})$ in this case.
Also, Bellingeri and Souli\'e \cite{Bellingeri-Soulie} defined this contruction for the welded braid group $WB_n$ which is a generalization of the braid group.
\end{remark}

\subsection{\textit{Reduced version}}
\begin{figure}[tbp]
\begin{center}
\includegraphics[height=100pt]{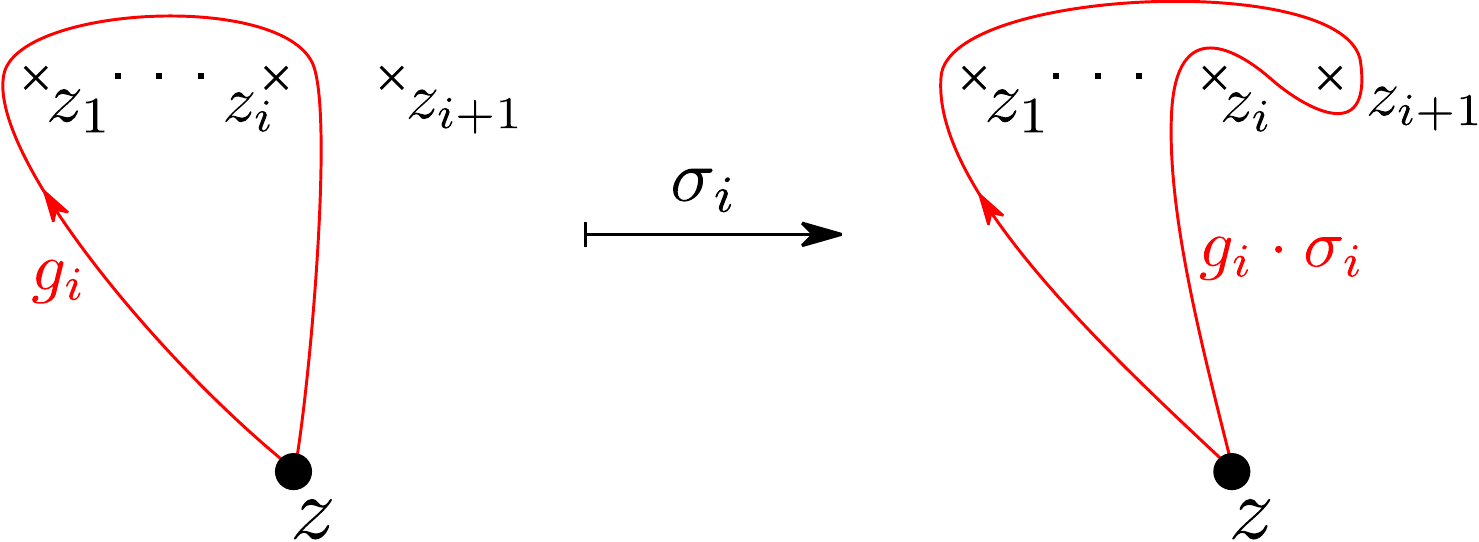}
\caption{The generator $g_i$ and the action of $\sigma_i$}
\label{braid_action_g}
\end{center}
\end{figure}

In order to describe a relation with the twisted Alexander invariant, we also consider other generators $g_1, \ldots, g_n$ of $F_n$, where $g_i := x_1 \cdots x_i$; see Figure \ref{braid_action_g}.
The action of $B_n$ on this new generating set is given by
\begin{equation*}
g_j \cdot \sigma_i :=
\begin{cases}
g_j & (j \neq i) \\
g_{i+1}g_i^{-1}g_{i-1} & (j = i \neq 1) \\
g_2 g_1^{-1} & (j = i = 1)
\end{cases}.
\end{equation*}
Since the augmentation ideal $\cI_{F_n}$ is also generated by $g_i - 1\ (1 \leq i \leq n)$, we have another matrix presentation of $\LM_{\ft_c} (\rho) (b)$ with respect to this generating set.
As mentioned in Section \ref{colored},  the loop $g_n = x_1 \cdots x_n$ is fixed by the action of $B_n$, and thus the matrix is written as
\begin{align*}
&\LM_{\ft_c} (\rho) (b) \\
&= \Diag \left( \rho (b), \ldots, \rho (b) \right) \cdot
\left(
\begin{array}{cc}
\left( (\rho \otimes \ft_c) \left( \dfr{\partial (g_i \cdot b)}{\partial g_j} \right) \right)_{1 \leq i,j \leq n-1} & V \\
0 & (\rho \otimes \ft_c) \left( \dfr{\partial (g_n \cdot b)}{\partial g_n} \right)
\end{array}
\right)\\
&= \Diag \left( \rho (b), \ldots, \rho (b) \right) \cdot
\left(
\begin{array}{cc}
\left( (\rho \otimes \ft_c) \left( \dfr{\partial (g_i \cdot b)}{\partial g_j} \right) \right)_{1 \leq i,j \leq n-1} & V \\
0 & I_k
\end{array}
\right),
\end{align*}
where $V$ is some $(n-1) k \times k$ matrix.
Therefore we obtain a subrepresentation of $\LM_{\ft_c} (\rho)$ by considering the first $n-1$ blocks:

\begin{defi}
Regard $F_{n-1}$ as being generated by $g_1, \ldots, g_{n-1}$.
The \textbf{reduced multivariable Long--Moody construction} of  the representation $\rho$ is a representation $\wtil{\LM}_{\ft_c} (\rho) \colon B_n \to GL_{R[H]}(R[H]^{\oplus k} \otimes_{R[F_{n-1}]} \cI_{F_{n-1}}) \cong GL_{(n-1)k}(R[H])$ given by
\begin{equation*}
\wtil{\LM}_{\ft_c} (\rho) (b) (v \otimes h) := \bar{\ft}_c (b) \cdot \left( v \left((\rho \otimes \ft_c) (b) \right) \otimes h \cdot b \right)
\end{equation*}
for any $b \in B_n$, $h \in \cI_{F_{n-1}}$ and $v \in R[H]^{\oplus k}$.
\end{defi}

Its matrix presentation with respect to the basis $\left\{ g_i -1 \relmiddle| 1 \leq i \leq n-1 \right\}$ is
\begin{equation*}
\wtil{\LM}_{\ft_c} (\rho) (b) = \Diag \left( \rho (b), \ldots, \rho (b) \right) \cdot \left( (\rho \otimes \ft_c) \left( \frac{\partial (g_i \cdot b)}{\partial g_j} \right) \right)_{1 \leq i,j \leq n-1}.
\end{equation*}

\begin{ex}
For the one-dimensional trivial representation $\cT \colon B_c \ltimes F_n \to GL_1(\bZ)$, the multivariable reduced Long--Moody construction $\wtil{\LM}_{\ft_c}(\cT) \colon B_c \to GL_{n-1}(\bZ[t_1^{\pm1}, \ldots, t_{\mu}^{\pm1}])$ is equivalent to the reduced Burau representation if $c = (1, 1, \ldots, 1)$, and equivalent to the reduced Gassner representation if $c = (1, 2, \ldots, n)$.
\end{ex}

\begin{remark} \label{twisted_homology}
The Long--Moody construction has an interpretation as the twisted homology.
Fix a sequence $c = (c_1, \ldots, c_n)$ of integers $1, 2, \ldots, \mu$.
Consider the surjection $\varphi_c \colon \pi_1(D_n, z) \to H$ given by $\varphi_c(x_i) = t_{c_i}$ for $1 \leq i \leq n$.
Note that $\varphi_c$ is equal to $\ft_c|_{F_n}$.
Let $p \colon \wtil{D}_n \to D_n$ be the universal covering of $D_n$ and set $\wtil{z} := p^{-1}(z)$.
Let $\rho \colon B_c \ltimes \pi_1(D_n, z) \left(\cong B_c \ltimes F_n \right) \to GL_k(R)$ be a representation.
By identifying $R^{\oplus k} \otimes_R R[H] \cong R[H]^{\oplus k}$ and restricting $\rho$ to $\pi_1(D_n, z)$, we obtain the representation $\rho \otimes \varphi_c \colon \pi_1(D_n, z) \to GL_k(R[H])$.
Then the chain complex of $R[H]$-modules
\begin{equation*}
C_*^{\rho \otimes \varphi_c} (D_n, z; R[H]^{\oplus k}) := R[H]^{\oplus k} \otimes_{\bZ[\pi_1(D_n, z)]} C_*(\wtil{D}_n, \wtil{z})
\end{equation*}
and the twisted homology group $H_*^{\rho \otimes \varphi_c} (D_n, z; R[H]^{\oplus k})$ are defined.
From \cite[Lemma 2.2]{Conway}, the $R[H]$-module $H_1^{\rho \otimes \varphi_c} (D_n, z; R[H]^{\oplus k})$ is free of rank $nk$.
As we have a representation $\rho|_{B_c} \colon B_c \to GL_k(R)$, the action $B_c \to Aut_{R[H]} (H_1^{\rho \otimes \varphi_c} (D_n, z; R[H]^{\oplus k})) \cong GL_{nk} (R[H])$ is well-defined, and it is equivalent to the multivariable Long--Moody construction $\LM_{\ft_c} (\rho) \colon B_c \to GL_{nk}(R[H])$.

On the other hand, Conway \cite{Conway} defined the twisted Burau map as a homomorphism on the twisted homology.
Namely, let $\rho \colon \pi_1 (D_n, z) \to GL_k(R)$ be a representation and $b$ a $(c, c)$-colored braid, then the braid $b$ induces a well-defined homomorphism
\begin{equation*}
\sB_{\rho} (b) \colon H_1^{b_* (\rho \otimes \varphi_c)} (D_n, z; R[H]^{\oplus k}) \to H_1^{\rho \otimes \varphi_c} (D_n, z; R[H]^{\oplus k}),
\end{equation*}
where $b_*$ is the isomorphism on the group $\pi_1(D_n, z)$ induced from $b$.
From the above isomorphism of the twisted homology, we obtain the map $\sB_{\rho} \colon B_c \to GL_{nk}(R[H])$, which is generally not a representation.
Its matrix presentation is of the form
\begin{equation*}
\sB_{\rho}(b) = \left( (\rho \otimes \varphi_c) \left( \frac{\partial (x_i \cdot b)}{\partial x_j} \right) \right)_{1 \leq i,j \leq n}.
\end{equation*}

In the setting above, we are able to consider the twisted homology group $H_*^{\rho \otimes \varphi_c} (D_n; R[H]^{\oplus k})$ and an action of $B_n$ on $H_1^{\rho \otimes \varphi_c} (D_n; R[H]^{\oplus k})$.
Conway \cite[Section 3.2]{Conway} studied this group in detail.
However, there is no obvious basis so that we are able to compute a matrix of $\wtil{\LM}_{\ft_c} (\rho)$.
\end{remark}

\subsection{\textit{Relation with the twisted Alexander invariant}}
Conway \cite[Theorem 3.15]{Conway} proved a relation with the reduced twisted Burau map and the twisted Alexander invariant.
The following theorem shows a relation with the multivariable reduced Long--Moody construction and the twisted Alexander invariant.
\begin{thm} \label{main}
Let $b \in B_c$ be a $\mu$-colored braid and $\hat{b}$ its closure.
Let $\rho \colon B_c \ltimes F_n \to GL_{k}(R)$ be a representation such that the restriction $\rho|_{F_n} \colon F_n \to GL_{k}(R)$ factors through the homomorphism $\phi$ and a representation $\pi_1 (S^3 \setminus \hat{b}) \to GL_{k}(R)$.
In other words, we assume the following diagram is commutative:
\begin{equation*}
\begin{tikzcd}
{F_n} && {GL_k (R)} \\
\\
{\pi_1 (S^3 \setminus \hat{b})}
\arrow["{\rho |_{F_n}}", from=1-1, to=1-3]
\arrow["\phi"', from=1-1, to=3-1]
\arrow[from=3-1, to=1-3]
\end{tikzcd}
\end{equation*}
Then
\begin{equation*}
\Delta_{\hat{b},\rho}(t_1, \ldots, t_{\mu}) = \pm \varepsilon t_{1}^{l_1} \cdots t_{\mu}^{l_{\mu}} \frac{\det \left( \wtil{\LM}_{\ft_c}(\rho) (b) - \Diag \left( \rho (b), \ldots, \rho (b) \right) \right)}{\det(\rho(x_1 \cdots x_n) t_{c_1} \cdots t_{c_n} -I_k)}
\end{equation*}
for some $\varepsilon \in R^{\times}$ and $l_1, \ldots, l_{\mu} \in \bZ$, where the representation $\pi_1 (S^3 \setminus \hat{b}) \to GL_{k}(R)$ is induced by $\rho|_{F_n}$, and let the same symbol $\rho$ denote this representation.
\end{thm}

\begin{proof}
Recall the presentation of $\pi_1 (S^3 \setminus \hat{b})$ using the braid action in Section \ref{colored}.
Replacing the generators $x_1, \ldots, x_n$ with $g_1, \ldots, g_n$, $\pi_1 (S^3 \setminus \hat{b})$ also has a similar presentation:
\begin{equation*}
\left\langle g_1 , \ldots , g_n \relmiddle| g_1 = g_1 \cdot b, \ldots, g_{n-1} = g_{n-1} \cdot b \right\rangle.
\end{equation*}
Then the surjective homomorphism $\alpha \colon \pi_1 (S^3 \setminus \hat{b}) \to H$ is given by $g_i \mapsto t_{c_1} \cdots t_{c_i}$ for $1 \leq i \leq n$.
Computing the Alexander matrix of this presentation, we have
\begin{equation*}
M = \left( \Phi \left( \frac{\partial (g_i \cdot b)}{\partial g_j} \right) - \delta_{ij} I_k \right)_{1 \leq i \leq n-1, 1 \leq j \leq n}.
\end{equation*}
Hence the twisted Alexander invariant is
\begin{equation*}
\Delta_{\hat{b},\rho}(t) = \frac{\det M_n}{\det \Phi (g_n - 1)} = \frac{\det \left( \Phi \left( \dfr{\partial (g_i \cdot b)}{\partial g_j} \right) - I_{(n-1)k} \right)}{\det \Phi (g_n - 1)}.
\end{equation*}
Note that since $\alpha(g_n) = t_{c_1} \cdots t_{c_n} \neq 1$, we have $\det \Phi (g_n - 1) \neq 0$; see Section \ref{twisted}.

On the other hand, since $\Phi = (\wtil{\rho} \otimes \wtil{\alpha}) \circ \wtil{\phi}$, $\ft_c|_{F_n} = \alpha \circ \phi$ and $\rho|_{F_n} = \rho \circ \phi$, we have
\begin{align*}
\wtil{\LM}_{\ft_c} (\rho) (b) &= \Diag \left( \rho (b), \ldots, \rho (b) \right) \cdot \left( (\rho \otimes \ft_c) \left( \frac{\partial (g_i \cdot b)}{\partial g_j} \right) \right) \\
&= \Diag \left( \rho (b), \ldots, \rho (b) \right) \cdot \left( \Phi \left( \frac{\partial (g_i \cdot b)}{\partial g_j} \right) \right).
\end{align*}
Therefore
\begin{align*}
&\det \left( \wtil{\LM}_{\ft_c} (\rho) (b) - \Diag \left( \rho (b), \ldots, \rho (b) \right) \right) \\
=& \det \Diag \left( \rho (b), \ldots, \rho (b) \right) \det \left( (\rho \otimes \ft_c) \left( \frac{\partial (g_i \cdot b)}{\partial g_j} \right) - I_{(n-1)k} \right) \\
=& (\det \rho (b))^{n-1} \det \left( \Phi \left( \frac{\partial (g_i \cdot b)}{\partial g_j} \right) - I_{(n-1)k} \right)
\end{align*}
and
\begin{equation*}
\det (\rho(x_1 \cdots x_n) t_{c_1} \cdots t_{c_n} -I_k) = \det ( (\rho \otimes \ft_c) (g_n) - I_k) = \det \Phi (g_n - 1).
\end{equation*}
Since $\det \rho (b) \in R^{\times}$, the theorem holds.
\end{proof}

\begin{remark}
In Remark \ref{twisted_homology}, if the input representation $\rho \colon \pi_1 (D_n, z) \to GL_k(R)$ extends to the semidirect product $B_c \ltimes \pi_1 (D_n, z)$, we can extend the twisted Burau map by using the representation $\rho|_{B_c} \colon B_c \to GL_k(R)$, that is, we can define the Long--Moody construction.
In this sense, the matrix from Theorem \ref{matrix} recovers the one of the twisted Burau map, and then Theorem \ref{main} could be a corollary of \cite[Theorem 1.1]{Conway}.
\end{remark}

\section{Examples} \label{ex}

In this section, we give some examples of Theorem \ref{main}.
The theorem supposes that a representation $\rho \colon B_c \ltimes F_n \to GL_{k}(R)$ satisfies the condition.
Conversely, let us choose a representation $\rho \colon \pi_1 (S^3 \setminus \hat{b}) \to GL_{k}(R)$ that extends to $B_c \ltimes F_n$ in the following examples.

\begin{ex} (\cite[Section 4]{Wada})
Set $b := \sigma_1^3 \in B_2 = B_{(1, 1)}$. Its closure is the trefoil knot $3_1$; see Figure \ref{trefoil}.
\begin{figure}[tbp]
\begin{center}
\includegraphics[height=80pt]{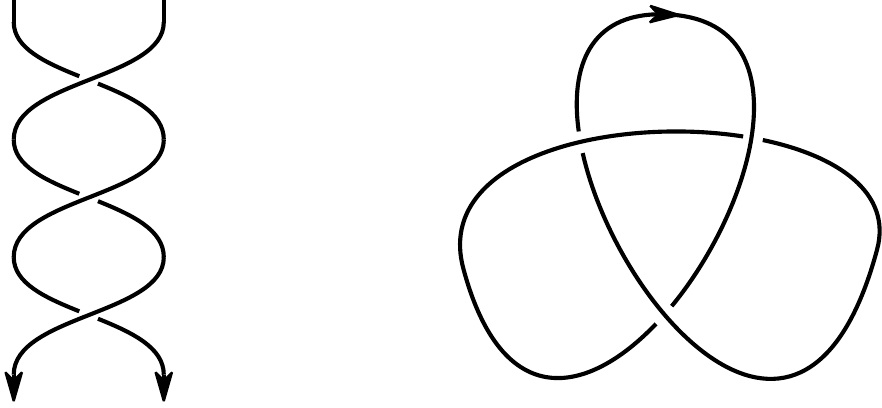}
\caption{The braid $\sigma_1^3$ and the trefoil knot $3_1$}
\label{trefoil}
\end{center}
\end{figure}
Then
\begin{align*}
\pi_1 (S^3 \setminus \hat{b})
&\cong \left\langle x_1, x_2 \relmiddle| x_1 = x_1 \cdot b, x_2 = x_2 \cdot b \right\rangle \\
&\cong \left\langle x_1, x_2 \relmiddle| x_1 x_2 x_1 = x_2 x_1 x_2 \right\rangle \\
&\cong \left\langle g_1, g_2 \relmiddle| g_1 =g_1 \cdot b \right\rangle.
\end{align*}
Let $\rho := \wtil{\mathscr{B}} \colon \pi_1 (S^3 \setminus \hat{b}) (\cong B_3) \to GL_2(\bZ[s^{\pm1}])$ be the reduced Burau representation given by
\begin{equation*}
\rho (x_1) =
\left(
\begin{array}{cc}
-s & 1 \\
0 & 1
\end{array}
\right)\ \ \textrm{and}\ \ 
\rho (x_2) =
\left(
\begin{array}{cc}
1 & 0 \\
s & -s
\end{array}
\right).
\end{equation*}
This is naturally extended to $F_2$.
The group $B_2 \ltimes F_2$ has the following presentation:
\begin{equation*}
\left\langle x_1, x_2, \sigma_1 \relmiddle| x_1 \sigma_1 = \sigma_1 x_1 x_2 x_1^{-1}, x_2 \sigma_1 = \sigma_1 x_1 \right\rangle,
\end{equation*}
and then if we set
\begin{equation*}
\rho (\sigma_1) :=
\left(
\begin{array}{cc}
0 & 1 \\
-1 & 1
\end{array}
\right),
\end{equation*}
the representation $\rho \colon B_2 \ltimes F_2 \to GL_{2}(\bZ[s^{\pm1}])$ is well-defined.
Since $g_1 \cdot b = g_2^2 g_1^{-1} g_2^{-1}$, we have
\begin{equation*}
\frac{\partial (g_1 \cdot b)}{\partial g_1} = -g_2^2 g_1^{-1}.
\end{equation*}
Therefore
\begin{equation*}
\wtil{\LM}_{\ft_c} (\rho) (b) = -\rho (b) \rho (g_2^2 g_1^{-1}) t^3 =
\left(
\begin{array}{cc}
s t^3 & - s t^3 + s^2 t^3 \\
s t^3 & - s t^3 \\
\end{array}
\right),
\end{equation*}
and thus we obtain
\begin{equation*}
\det \left( \wtil{\LM}_{\ft_c} (\rho) (b) - \rho(b) \right) =
\det \left(
\begin{array}{cc}
1 + s t^3 & - s t^3 + s^2 t^3 \\
s t^3 & 1 - s t^3 \\
\end{array}
\right)
= 1 - s^3 t^6
\end{equation*}
and
\begin{equation*}
\det \left( \rho(x_1 x_2) t^2 - I_2 \right) =
\det \left(
\begin{array}{cc}
-1 & -s t^2 \\
s t^2 & -s t^2-1 \\
\end{array}
\right)
= 1 + s t^2 + s^2 t^4.
\end{equation*}
Hence
\begin{equation*}
\Delta_{\hat{b},\rho}(t) = \frac{1 - s^3 t^6}{1 + s t^2 + s^2 t^4} = 1 -s t^2.
\end{equation*}
\end{ex}

\begin{ex} \label{figure-8}
Set $b := \sigma_1 \sigma_2^{-1} \sigma_1 \sigma_2^{-1} \in B_3 = B_{(1, 1, 1)}$. Its closure is the figure eight knot $4_1$; see Figure \ref{figure_eight}.
\begin{figure}[tbp]
\begin{center}
\includegraphics[height=100pt]{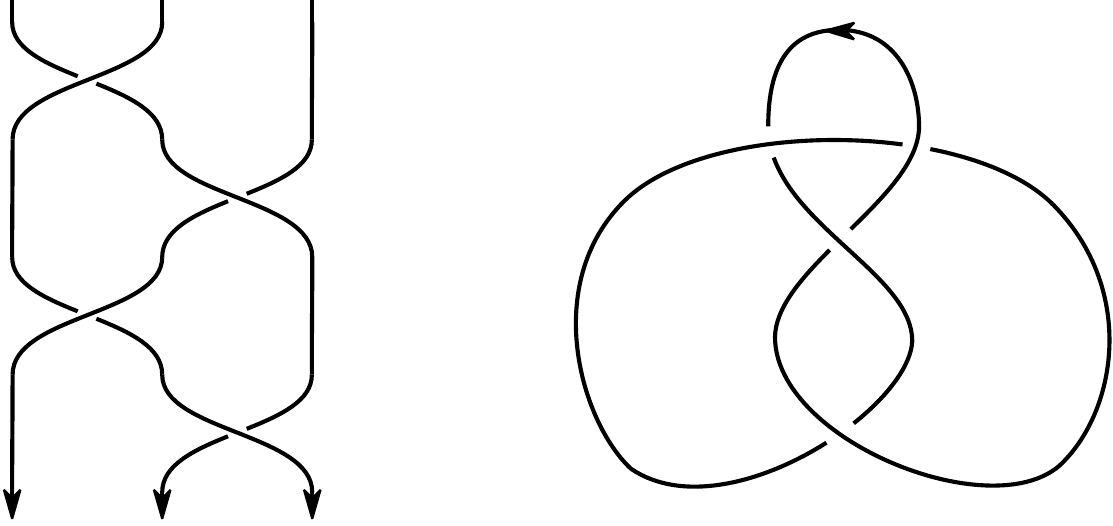}
\caption{The braid $\sigma_1 \sigma_2^{-1} \sigma_1 \sigma_2^{-1}$ and the figure eight knot $4_1$}
\label{figure_eight}
\end{center}
\end{figure}
Then
\begin{align*}
\pi_1 (S^3 \setminus \hat{b})
&\cong \left\langle x_1, x_2, x_3 \relmiddle| x_1 = x_1 \cdot b, x_2 = x_2 \cdot b, x_3 = x_3 \cdot b \right\rangle \\
&\cong \left\langle x_2, x_3 \relmiddle| x_2 [x_3^{-1}, x_2] = [x_3^{-1}, x_2] x_3 \right\rangle \\
&\cong \left\langle g_1, g_2, g_3 \relmiddle| g_1 = g_1 \cdot b, g_2 = g_2 \cdot b \right\rangle,
\end{align*}
where $[x, y] := x y x^{-1} y^{-1}$.
Let $\rho \colon \pi_1 (S^3 \setminus \hat{b}) \to SL_2(\bZ_7)$ be the representation given by
\begin{equation*}
\rho (x_1) =
\left(
\begin{array}{cc}
2 & 0 \\
0 & 4
\end{array}
\right),\ \ 
\rho (x_2) =
\left(
\begin{array}{cc}
1 & 1 \\
4 & 5
\end{array}
\right)\ \ \textrm{and}\ \ 
\rho (x_3) =
\left(
\begin{array}{cc}
1 & 2 \\
2 & 5
\end{array}
\right).
\end{equation*}
The group $B_3 \ltimes F_3$ has a presentation
\begin{equation*}
\left\langle x_1, x_2, x_3, \sigma_1, \sigma_2 \relmiddle|
\begin{array}{l}
x_1 \sigma_1 = \sigma_1 x_1 x_2 x_1^{-1}, x_2 \sigma_1 = \sigma_1 x_1, x_3 \sigma_1 = \sigma_1 x_3, \\
x_1 \sigma_2 = \sigma_2 x_1, x_2 \sigma_2 = \sigma_2 x_2 x_3 x_2^{-1}, x_3 \sigma_2 = \sigma_2 x_2, \\
\sigma_1 \sigma_2 \sigma_1 = \sigma_2 \sigma_1 \sigma_2
\end{array}
\right\rangle,
\end{equation*}
and setting
\begin{equation*}
\rho (\sigma_1) =
\left(
\begin{array}{cc}
5 & 5 \\
5 & 1
\end{array}
\right)\ \ \textrm{and}\ \ 
\rho (\sigma_2) =
\left(
\begin{array}{cc}
4 & 0 \\
0 & 2
\end{array}
\right),
\end{equation*}
the representation $\rho$ is extended to $B_3 \ltimes F_3$.
From $g_1 \cdot b = g_1 g_2^{-1} g_3 g_1^{-1} g_3^{-1} g_2 g_1^{-1} g_3 g_1 g_3^{-1} g_2 g_1^{-1}$ and $g_2 \cdot b = g_1 g_2^{-1} g_3 g_1^{-1} g_3^{-1} g_2 g_1^{-1} g_3$,
\begin{align*}
\begin{split}
\frac{\partial (g_1 \cdot b)}{\partial g_1} &= 1 - g_1 g_2^{-1} g_3 g_1^{-1} - g_1 g_2^{-1} g_3 g_1^{-1} g_3^{-1} g_2 g_1^{-1}\\
&\quad \quad + g_1 g_2^{-1} g_3 g_1^{-1} g_3^{-1} g_2 g_1^{-1} g_3 - g_1 g_2^{-1} g_3 g_1^{-1} g_3^{-1} g_2 g_1^{-1} g_3 g_1 g_3^{-1} g_2 g_1^{-1},\\
\frac{\partial (g_1 \cdot b)}{\partial g_2} &= - g_1 g_2^{-1} + g_1 g_2^{-1} g_3 g_1^{-1} g_3^{-1} + g_1 g_2^{-1} g_3 g_1^{-1} g_3^{-1} g_2 g_1^{-1} g_3 g_1 g_3^{-1},\\
\frac{\partial (g_2 \cdot b)}{\partial g_1} &= 1 - g_1 g_2^{-1} g_3 g_1^{-1} - g_1 g_2^{-1} g_3 g_1^{-1} g_3^{-1} g_2 g_1^{-1},\\
\frac{\partial (g_2 \cdot b)}{\partial g_2} &= - g_1 g_2^{-1} + g_1 g_2^{-1} g_3 g_1^{-1} g_3^{-1}.
\end{split}
\end{align*}
Therefore
\begin{align*}
\wtil{\LM}_{\ft_c} (\rho) (b) &= 
\left(
\begin{array}{cc}
\rho (b) & 0 \\
0 & \rho (b) \\
\end{array}
\right)
\left(
\begin{array}{cc}
(\rho \otimes \ft_c) \left( \dfr{\partial (g_1 \cdot b)}{\partial g_1} \right) & (\rho \otimes \ft_c) \left( \dfr{\partial (g_1 \cdot b)}{\partial g_2} \right) \\
(\rho \otimes \ft_c) \left( \dfr{\partial (g_2 \cdot b)}{\partial g_1} \right) & (\rho \otimes \ft_c) \left( \dfr{\partial (g_2 \cdot b)}{\partial g_2} \right) \\
\end{array}
\right)\\
&= \left( \large
\begin{array}{cccc}
\frac{5t^3+3t^2+6t+5}{t} & \frac{5t^3+t^2+5}{t} & \frac{4t^2+5t+1}{t^2} & \frac{t^2+3 t+1}{t^2} \\
\frac{5t^3+4t^2+5}{t} & \frac{t^3+2t^2+6t+1}{t} & \frac{4t^2+6t+4}{t^2} & \frac{3t^2+t+5}{t^2} \\
\frac{t^2+6t+5}{t} & \frac{t^2+5}{t} & \frac{5t+1}{t^2} & \frac{3t+1}{t^2} \\
\frac{4t^2+5}{t} & \frac{5t^2+6 t+1}{t} & \frac{6t+4}{t^2} & \frac
{t+5}{t^2} \\
\end{array}
\right).
\end{align*}
Hence
\begin{align*}
&\det \left( \wtil{\LM}_{\ft_c} (\rho) (b) -
\left(
\begin{array}{cc}
\rho (b) & 0 \\
0 & \rho (b) \\
\end{array}
\right)
\right) \\
= &\det \left( \large
\begin{array}{cccc}
\frac{5t^3+3t^2+5}{t} & \frac{5t^3+t^2+5}{t} & \frac{4t^2+5t+1}{t^2} & \frac{t^2+3t+1}{t^2} \\
\frac{5 t^3+4 t^2+5}{t} & \frac{t^3+2t^2+1}{t} & \frac{4t^2+6 t+4}{t^2} & \frac{3t^2+t+5}{t^2} \\
\frac{t^2+6t+5}{t} & \frac{t^2+5}{t} & \frac{t^2+5t+1}{t^2} & \frac{3t+1}{t^2} \\
\frac{4t^2+5}{t} & \frac{5t^2+6t+1}{t} & \frac{6t+4}{t^2} & \frac{t^2+t+5}{t^2} \\
\end{array}
\right) \\
= &(t + 1)^4 (t + 2)^2 (t + 4)^2 t^{-4}
\end{align*}
and
\begin{equation*}
\det \left( \rho(x_1 x_2 x_3) t^3 - I_2 \right) =
\det \left(
\begin{array}{cc}
6 + 6t^3 & 0 \\
0 & 6 + 6t^3 \\
\end{array}
\right)
= (t + 1)^2 (t + 2)^2 (t + 4)^2.
\end{equation*}
Consequently
\begin{equation*}
\Delta_{\hat{b},\rho}(t) = (t + 1)^2.
\end{equation*}
\end{ex}

\begin{ex} \label{figure-8_2}
Consider the same braid in Example \ref{figure-8}.
Let $\rho \colon \pi_1 (S^3 \setminus \hat{b}) \to SL_2(\bC)$ be the representation given by
\begin{equation*}
\rho (x_1) =
\dfr{1}{3} \left(
\begin{array}{cc}
-(s+2) & 3s \\
-2s^2 & s-1
\end{array}
\right),\ \ 
\rho (x_2) =
\left(
\begin{array}{cc}
s & 0 \\
0 & s^{-1}
\end{array}
\right)\ \ \textrm{and}\ \ 
\rho (x_3) =
\dfr{1}{3} \left(
\begin{array}{cc}
-(s+2) & 3 \\
-2 & s-1
\end{array}
\right),
\end{equation*}
where $s \in \bC$ satisfies $s^2 + s + 1 = 0$.
Setting
\begin{equation*}
\rho (\sigma_1) =
\sqrt{\dfr{-s}{3}} \left(
\begin{array}{cc}
1 & s+2 \\
-\displaystyle{\frac{2(s+2)}{3}} & -s^2
\end{array}
\right)\ \ \textrm{and}\ \ 
\rho (\sigma_2) =
\sqrt{\dfr{-s}{3}} \left(
\begin{array}{cc}
1 & s-1 \\
\displaystyle{\frac{2(2s+1)}{3}} & -s^2
\end{array}
\right),
\end{equation*}
the representation $\rho$ is extended to $B_3 \ltimes F_3$.
Then
\begin{align*}
&\det \left( \wtil{\LM}_{\ft_c} (\rho) (b) -
\left(
\begin{array}{cc}
\rho (b) & 0 \\
0 & \rho (b) \\
\end{array}
\right)
\right) \\
=& \det \left( \large
\begin{array}{cccc}
\frac{\left(t^3+2t^2+1\right) (s-1)}{3t} & -\frac{t^3-t^2s+1}{t} & \frac{2t^2s+t^2+ts-t+3s}{3t^2} & \frac{ts+s+1}{t} \\
\frac{2 \left(-t^2s^2+t^3+1\right)}{3t} & -\frac{\left(t^3+2t^2+1\right) (s+2)}{3t} & \frac{2 (ts+t+s)}{3t} & -\frac{2t^2s+t^2-3s^2+ts+2t}{3t^2} \\
\frac{3t^2s-3t+s-1}{3t} & -\frac{1}{t} & \frac{3t^2+ts-t+3s}{3t^2} & \frac{s+1}{t} \\
\frac{2}{3t} & \frac{3t^2s^2-3t-s-2}{3t} & \frac{2s}{3t} & \frac{3t^2+3s^2-ts-2t}{3t^2}\\
\end{array}
\right)\\
=& \dfr{(1+t)^4 (1-t+t^2)^2}{t^4},
\end{align*}
and
\begin{align*}
\det \left( \rho(x_1 x_2 x_3) t^3 - I_2 \right) &=
\det \left(
\begin{array}{cc}
-(1+t) (1-t+t^2) & 0 \\
0 & -(1+t) (1-t+t^2) \\
\end{array}
\right)\\
&= (1+t)^2 (1-t+t^2)^2.
\end{align*}
Therefore
\begin{equation*}
\Delta_{\hat{b},\rho}(t) = (t + 1)^2.
\end{equation*}
\end{ex}

\begin{remark}
Let $\rho \colon \pi_1 (S^3 \setminus 4_1) \to SL_2(\bC)$ be the holonomy representation of the figure eight knot $4_1$ given by
\begin{equation*}
\rho (x_2) =
\left(
\begin{array}{cc}
1 & 1 \\
0 & 1
\end{array}
\right)\ \ \textrm{and}\ \ 
\rho (x_3) =
\left(
\begin{array}{cc}
1 & 0 \\
-s & 1
\end{array}
\right),
\end{equation*}
where $s \in \bC$ satisfies $s^2 + s + 1 = 0$.
However, this representation does not extend to $B_3 \ltimes F_3$.
\end{remark}

\begin{remark}
In general, there is a canonical inclusion $\iota \colon B_n \ltimes F_n \hookrightarrow B_{n+1}$ with the correspondence $\sigma_i \mapsto \sigma_{i+1}$ and $x_i \mapsto \sigma_1^{-1} \cdots \sigma_{i-1}^{-1} \sigma_i^2 \sigma_{i-1} \cdots \sigma_1$.
Under some mild assumption, a representation $\rho \colon B_n \ltimes F_n \to GL_k(\bC)$ can be extended to a representation $\rho' \colon B_{n+1} \to GL_{k}(\bC)$ as follows.
First, we have $\rho'(\sigma_{i+1}) := \rho(\sigma_i)$ for $1 \leq i \leq n-1$.
Then, defining an appropriate image for the generator $\sigma_1$ is enough to define $\rho'$.
Note that $\iota(x_1) = \sigma_1^2$ and every element of $GL_k(\bC)$ has a square root; see for example \cite{Bjorck-Hammarling, Cross-Lancaster}.
Hence we assume that there exists a square root $A$ of $\rho(x_1)$ which satisfies the equality $A \rho'(\sigma_2) A = \rho'(\sigma_2) A \rho'(\sigma_2)$, and define $\rho'(\sigma_1)$ as the matrix $A$.

For example, the representation $\rho \colon B_3 \ltimes F_3 \to SL_2(\bC)$ given in Example \ref{figure-8_2} extends to $B_4$ and is essentially equivalent to $\wtil{\mathscr{B}} \circ \pi \colon B_4 \to SL_2(\bC)$, where $\pi \colon B_4 \to B_3$ is given by $\sigma_1, \sigma_3 \mapsto \sigma_1$ and $\sigma_2 \mapsto \sigma_2$, and the reduced Burau representation $\wtil{\mathscr{B}}$ of $B_3$ is regarded as a complex representation by taking the tensor product $\mathbb{C} \otimes \mathbb{Z} [s^{\pm 1}] \cong \mathbb{C} [s^{\pm 1}]$ and specializing $s$ to a non-zero complex value.
\end{remark}

\begin{ex} (\cite[Theorem 4.2]{Kitano-Morifuji1})
Consider the 2-braid $b := \sigma_1^q \in B_2 = B_{(1, 1)}$, where $q = 2m + 1$ $(m \in \bZ_{> 0})$.
Its closure is the $(2,q)$-torus knot $T(2,q)$.
Then
\begin{align*}
\pi_1 (S^3 \setminus \hat{b})
&\cong \left\langle x, y \relmiddle| x^2 = y^q \right\rangle \\
&\cong \left\langle g_1, g_2 \relmiddle| g_1 = g_1 \cdot b = g_2^{m+1} g_1^{-1} g_2^{-m} \right\rangle \\
&\cong \left\langle x_1, x_2 \relmiddle| x_1 = x_1 \cdot b = (x_1 x_2)^{m+1} x_1^{-1} (x_1 x_2)^{-m}, x_2 = x_2 \cdot b = (x_1 x_2)^m x_1 (x_1 x_2)^{-m} \right\rangle \\
&\cong \left\langle x_1, x_2 \relmiddle| (x_1 x_2)^m x_1 = x_2 (x_1 x_2)^m \right\rangle,
\end{align*}
where $x$ represents the core of the solid torus such that the torus knot $T(2,q)$ is stuck on its boundary torus $T \subset S^3$ and $y$ represents that of another solid torus whose boundary is $T$, namely $T$ defines the genus one Heegaard decomposition of $S^3$ by two solid tori; see Figure \ref{(2,q)_torus}.
The generators $x$ and $y$ are written in terms of $g_1$ and $g_2$ as follows, respectively:

\begin{figure}[tbp]
\begin{center}
\includegraphics[height=120pt]{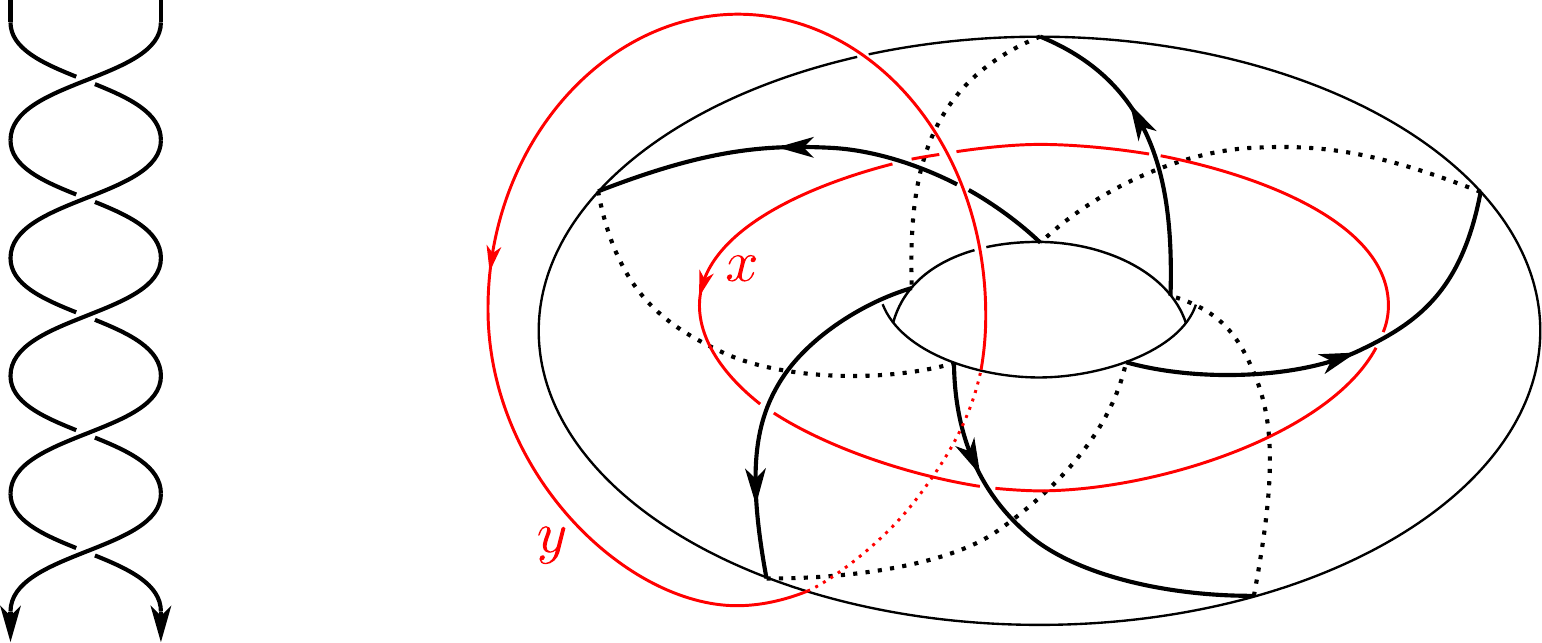}
\caption{The braid $\sigma_1^q$ and the $(2,q)$-torus knot $T(2,q)$ ($q=5$)}
\label{(2,q)_torus}
\end{center}
\end{figure}

\begin{align} \label{rel}
\begin{cases}
x = x_2 (x_2 \cdot \sigma_1^{-1}) \cdots (x_2 \cdot \sigma_1^{-(q-1)}) = x_2 (x_1 x_2)^m = g_1^{-1} g_2^{m+1}, \\
y = x_1 x_2 = g_2.
\end{cases}
\end{align}
Consider the non-abelian irreducible representation $\rho \colon \pi_1 (S^3 \setminus \hat{b}) \to SL_2(\bC)$ given by
\begin{equation*}
\rho (x) =
\left(
\begin{array}{cc}
s & 1 \\
u & v
\end{array}
\right)\ \ \textrm{and}\ \ 
\rho (y) =
\left(
\begin{array}{cc}
\xi_r & 0 \\
0 & \xi^{-1}_r
\end{array}
\right),
\end{equation*}
where $u = s v - 1 \neq 0$, $r$ is odd $(0 < r < q)$ and $\xi_r = \exp \left(\dfr{\sqrt{-1} \pi r}{q} \right)$.
By the relation $x^2 = y^q$, it holds that $\rho (x)^2 = \rho (y)^q = -I_2$ and we obtain $u = - 1 - s^2$ and $v = -s$.
From (\ref{rel}), we have
\begin{align*}
\rho (g_1) &= \rho (y)^{m+1} \rho (x)^{-1} =
\left(
\begin{array}{cc}
-\sqrt{-1}^r \xi^{\frac{1}{2}}_{r} s & -\sqrt{-1}^r \xi^{\frac{1}{2}}_{r} \\
\sqrt{-1}^{-r} \xi^{-\frac{1}{2}}_{r} (1 + s^2) & \sqrt{-1}^{-r} \xi^{-\frac{1}{2}}_{r} s
\end{array}
\right), \\
\rho (g_2) &= \rho (y) =
\left(
\begin{array}{cc}
\xi_r & 0 \\
0 & \xi^{-1}_r
\end{array}
\right).
\end{align*}
We regard $F_2$ as being generated by $g_1$ and $g_2$, and then a presentation of $B_2 \ltimes F_2$ is written as
\begin{equation*}
\left\langle g_1, g_2, \sigma_1 \relmiddle| g_1 \sigma_1 = \sigma_1 g_2 g_1^{-1}, g_2 \sigma_1 = \sigma_1 g_2 \right\rangle.
\end{equation*}
If we set
\begin{equation*}
\rho (\sigma_1) :=
\left(
\begin{array}{cc}
\sqrt{-1} \xi^{-\frac{1}{2}}_{r} & 0 \\
0 & -\sqrt{-1} \xi^{\frac{1}{2}}_{r}
\end{array}
\right),
\end{equation*}
this representation is extended to $B_2 \ltimes F_2$.
Since $g_1 \cdot b = g_2^{m+1} g_1^{-1} g_2^{-m}$, we have
\begin{equation*}
\frac{\partial (g_1 \cdot b)}{\partial g_1} = -g_2^{m+1} g_1^{-1},
\end{equation*}
and thus
\begin{align*}
\wtil{\LM}_{\ft_c} (\rho) (b) &= -\rho (\sigma_1^q) \rho (g_2^{m+1} g_1^{-1}) t^{2m+1} \\
&= -\left(
\begin{array}{cc}
\sqrt{-1}^{q-r} & 0 \\
0 & -\sqrt{-1}^{q+r} \\
\end{array}
\right)
\left(
\begin{array}{cc}
s & -\xi_{r} \\
\xi^{-1}_{r} (1 + s^2) & -s
\end{array}
\right) t^{2m+1} \\
&= \left(
\begin{array}{cc}
-\sqrt{-1}^{q-r} s t^{2m+1} & \sqrt{-1}^{q-r} \xi_{r} t^{2m+1} \\
\sqrt{-1}^{q+r} \xi^{-1}_{r} (1 + s^2) t^{2m+1} & -\sqrt{-1}^{q+r} s t^{2m+1}
\end{array}
\right).
\end{align*}
Therefore
\begin{align*}
\det \left( \wtil{\LM}_{\ft_c} (\rho) (b) - \rho(b) \right) &=
\det \left(
\begin{array}{cc}
-\sqrt{-1}^{q-r} (1 + s t^{2m+1}) & \sqrt{-1}^{q-r} \xi_{r} t^{2m+1} \\
\sqrt{-1}^{q+r} \xi^{-1}_{r} (1 + s^2) t^{2m+1} & \sqrt{-1}^{q+r} (1 - s t^{2m+1})
\end{array}
\right) \\
&= -\sqrt{-1}^{2q} - \sqrt{-1}^{2q} t^{2m+4} = 1 + t^{2q}
\end{align*}
and
\begin{equation*}
\det \left( \rho(x_1 x_2) t^2 - I_2 \right) =
\det \left(
\begin{array}{cc}
-1 + \xi_{r} t^2 & 0 \\
0 & -1 + \xi^{-1}_{r} t^2 \\
\end{array}
\right)
= \left( 1 - \xi_{r} t^2 \right) \left( 1 - \xi^{-1}_{r} t^2 \right).
\end{equation*}
Hence
\begin{equation*}
\Delta_{\hat{b},\rho}(t) = \frac{1 + t^{2q}}{\left( 1 - \xi_{r} t^2 \right) \left( 1 - \xi^{-1}_{r} t^2 \right)}
\end{equation*}
and this is a Laurent polynomial in $\bC[t^{\pm1}]$, justifying the result that the twisted Alexander invariant associated with a non-abelian $SL_2(\bC)$-representation is a Laurent polynomial by Kitano and Morifuji \cite[Theorem 3.1]{Kitano-Morifuji2}.
\end{ex}

\section*{Acknowledgments}
The author would like to thank Takuya Sakasai for his careful reading of the paper and his helpful advice on this research.
He also would like to thank Arthur Souli\'e for his proofreading this paper and valuable comments about the Long--Moody construction.


\bibliography{LM_and_twisted}
\bibliographystyle{plain}

\end{document}